\documentclass[11pt]{article}
\usepackage{graphicx}
\usepackage{amssymb,amsmath,amsfonts,amsthm,amscd,mathrsfs}
\usepackage{epstopdf}
\usepackage{amsbsy}
\usepackage{psfrag, color}
\input xy
\xyoption{all}
\CompileMatrices

\numberwithin{equation}{section}

\newtheorem{theorem}[equation]{Theorem}
\newtheorem{lemma}[equation]{Lemma}

\newtheorem{corollary}[equation]{Corollary}
\newtheorem{proposition}[equation]{Proposition}
\theoremstyle{definition}

\newtheorem{defn/notation}[equation]{Definition/Notation}

\newtheorem{remark}[equation]{Remark}

\newtheorem{remark/definition}[equation]{Remark-Definition}

\newcommand{\ct}{\ensuremath{c_{12}}}
\newcommand{\sigm}{\ensuremath{\sigma^{\overline{m}}}}
\newcommand{\sign}{\ensuremath{\sigma^{\overline{n}}}}
\newcommand{\Z}{\ensuremath{\mathbb{Z}}}
\newcommand{\N}{\ensuremath{\mathbb{N}}}
\newcommand{\ig}{\ensuremath{I[G]}}

\newcommand{\atg}{\ensuremath{A_2(G_{12})}}
\newcommand{\barm}{\ensuremath{\overline{m}}}
\newcommand{\barn}{\ensuremath{\overline{n}}}

\newcommand{\norm}{\ensuremath{\textrm{N}}}
\newcommand{\ctn}[1]{\ensuremath{c_{#1}}}
\newcommand{\mo}{\ensuremath{M_{\omega}}}
\newcommand{\amm}{\ensuremath{\atg \oplus K_3 \oplus K_4}}
\newcommand{\res}{\ensuremath{\mathrm{res}}}
\newcommand{\zg}{\ensuremath{\Z[G]}}
\newcommand{\zgot}{\ensuremath{\Z[G_{12}]}}

\newcommand{\eqnref}[1]{(\ref{#1})}
\newcommand{\gal}{\ensuremath{\textrm{Gal}}}
\newcommand{\slashfrac}[2]{\ensuremath{\raise1ex\hbox{#1}\kern-.2em/\kern-.30em\lower1ex\hbox{#2}}}

\newcommand{\lot}{\ensuremath{L_{12}^{G_{12}}}}
\newcommand{\p}{\ensuremath{|\res^G_H(\omega)|}}
\newcommand{\zp}{\ensuremath{\Z/p\Z}}
\newcommand{\chara}{\ensuremath{\textrm{char}}}
\newcommand{\ind}{\ensuremath{\mathrm{ind}}}
\newcommand{\G}{\ensuremath{\langle \sigma_1\rangle\times\ldots\times\langle\sigma_r\rangle}}
\newcommand{\F}{\ensuremath{{\mathcal{F}_A}}}
\renewcommand{\inf}{\ensuremath{\mathrm{inf}}}
\newcommand{\Br}{\ensuremath{\mathrm{Br}}}
\newcommand{\br}{\ensuremath{\mathrm{Br}}}
\newcommand{\dec}{\ensuremath{\mathrm{Dec}}}

\title{Degeneracy and Decomposability in Abelian Crossed Products}
\author{Kelly McKinnie}
\begin{document}
\maketitle

\abstract{In this paper we continue the study of the relationship between degeneracy and decomposability in abelian crossed products (\cite{mck-galoissubfields}).  In particular we construct an indecomposable abelian crossed product division algebra of exponent $p$ and index $p^2$ for $p$ an odd prime.  The algebra we construct is generic in the sense of \cite{AS} and has the property that its underlying abelian crossed product is a decomposable division algebra defined by a non-degenerate matrix.  This algebra gives an example of an indecomposable generic abelian crossed product which is shown to be indecomposable without using torsion in the Chow group of the corresponding Severi-Brauer variety as was needed in \cite{Karpenko2} and \cite{mck-galoissubfields}.  It also gives an example of a Brauer class which is in Tignol's Dec group with respect to one abelian maximal subfield, but not in the Dec group with respect to another.}

\section{Introduction}  
A finite dimensional division algebra $D$ with center a field $F$ is said to be {\it decomposable} if there exists an $F$-isomorphism $D \cong D_1\otimes_F D_2$ with $\ind(D_i)>1$.  By the primary decomposition theorem (see e.g., \cite[Cor.11, pg 68]{Draxl}), an indecomposable division algebra necessarily has prime power index.  Recall $\ind(D)=\sqrt{\dim_FD}$.  Moreover, if an $F$-division algebra $D$ has exponent equal to its index, then $D$ is indecomposable.  Therefore, it is only of interest to study indecomposable division algebras, $D$, of prime power index with $\ind(D)>\exp(D)$.  

Prime power index indecomposable division algebras have been studied in many contexts since a construction of one with $\ind(D)<\exp(D)$ was given by Saltman in \cite{Saltman-inde}.  For example, in light of the present paper, we draw the readers attention to \cite{mck-galoissubfields} and \cite{mounirh1} where the decomposability of generic abelian crossed product division algebras are studied.  In \cite[Cor. 3.6]{mounirh1} it is shown that a $p$-power degree generic abelian crossed product defined by a non-degenerate matrix (definition given below) is indecomposable.  \cite{mck-galoissubfields} and \cite{mounirh1} leave open the question of whether any abelian crossed product defined by a non-degenerate matrix is indecomposable.  In this paper we answer this question in the negative by constructing a decomposable abelian crossed product $\Delta_\F$  defined by a non-degenerate matrix (Theorem \ref{t1}/\ref{p1}).  As mentioned in more detail below, we show $\Delta_\F$ is defined by a non-degenerate matrix without using the Chow group as was done in \cite{Karpenko2} and \cite{mck-galoissubfields}.

Our example has two additional properties.  First, recall that in \cite[Prop. 3.1.1]{mck-galoissubfields} it is shown that an abelian crossed product $\Delta$ of exponent $p$ defined by a degenerate matrix has no nontrivial torsion in $\mathrm{CH}^2(X)$.  Here $X=SB(\Delta)$, the Severi-Brauer variety of $\Delta$.  Our abelian crossed product, $\Delta_\F$, provides a counterexample to the converse of this statement.  That is, our example satisfies $\mathrm{CH}^2(X)$ is torsion free where $X=SB(\Delta_\F)$ and $\Delta_\F$ is an exponent $p$ abelian crossed product defined by a non-degenerate matrix (see remark \ref{remark3} for more details).  As a consequence we see that the existence of non-trivial torsion in $\mathrm{CH}^2(X)$ for $X$ the Severi-Brauer variety of an abelian crossed product is not strictly controlled by degeneracy of the matrix of the abelian crossed product.

Secondly, our abelian crossed product $\Delta_\F$ with center $\F$ gives an example of a finite dimensional division algebra with two abelian maximal subfields, $N_1$ and $N_2$ such that $[\Delta_\F]\in \dec(N_1/\F)$ and $[\Delta_\F] \notin \dec(N_2/\F)$ (Corollaries \ref{cor2} and \ref{cor3}).  Here $\dec(\,\,)$ is the decomposition group of Tignol defined in \cite{Tignol-ACP}.  Recall for any field $F$ and finite abelian extension $N/F$, $\dec(N/F)$ is the subgroup of the relative Brauer group $\br(N/F)$ generated by the subgroups $\Br(K/F)$ where $K$ is a cyclic subfield of $N$.  Let $G=\gal(N/F)$.  Since $N/F$ is an abelian extension there is a basis $\{\sigma_i\}$ such that $G \cong \langle \sigma_1\rangle\times \cdots \times \langle \sigma_r\rangle$.  In \cite[Cor. 1.4]{Tignol-ACP} Tignol shows that $\dec(N/F)$ can also be described as the set of Brauer classes $[A]$ such that the algebra $A$ contains $N$, has $\deg(A)=[N:F]$ and such that $A$ decomposes into the tensor product of cyclic algebras
\[A \cong A_1\otimes_F\cdots\otimes_FA_r\]
 where for each $i$, the algebra $A_i$ is a cyclic $F$-algebra containing $K_i$ as a maximal subfield.  Here $K_i/F$ is a cyclic extension with $K_i=N^{G_i}$ and $G_i=\langle \sigma_1\rangle\times\cdots\times\widehat{\langle \sigma_i\rangle} \times \cdots\times\langle \sigma_r\rangle$.  If $[A] \in \dec(N/F)$, then $A$ is said to decompose ``according to $N/F$''.  The example in this paper shows that there exist Brauer classes $[D]$ such that $D$ decomposes with respect to one abelian maximal subfield, but not with respect to another. 

\subsection{Abelian Crossed Products and related definitions}
Let $F$ be a field.  An abelian crossed product is a central simple $F$-algebra which contains a maximal subfield that is abelian Galois over $F$.  Let $\Delta$ be an abelian crossed product over $F$ (we will write this as $\Delta/F$) with finite abelian maximal subfield $K$ and $G=\gal(K/F)=\G$.  As detailed in \cite{AS} or \cite{mck-primetop}, for every abelian crossed product there is a matrix $u=(u_{ij})\in M_r(K^*)$ and a vector $b=(b_i)_{i=1}^r\in(K^*)^r$ so that $\Delta$ is isomorphic to the following algebra.
\begin{equation}
\Delta \cong (K/F,G,z,u,b)=\bigoplus_{0\leq i_j \leq n_j}Kz_1^{i_1}\ldots z_r^{i_r}\label{eqn1}\end{equation}
Here $n_j=|\sigma_j|$ and multiplication in this algebra is given by the conditions $z_iz_j=u_{ij}z_jz_i$, $z_i^{n_i}=b_i$ and $z_ik=\sigma_i(k)z_i$ for all $k \in K$.  Throughout this paper we will use multi-index notation:  for $\barm=(m_1,\ldots,m_r) \in \N^r$ set $z^{\barm}=z_1^{m_1}\ldots z_r^{m_r}$ and $\sigma^{\barm}=\sigma_1^{m_1}\ldots \sigma_r^{m_r}$.  Moreover, set $u_{\barm,\barn}=z^{\barm}z^{\barn}(z^{\barm})^{-1}(z^{\barn})^{-1}\in K^*$.

Since the matrix $u$ determines multiplication in the algebra, properties of the matrix $u$ determine some properties of the abelian crossed product $\Delta$.  In \cite{AS} the notion of a matrix being degenerate was defined and this notion was further studied and extended in \cite{mck-primetop} and \cite{mounirh1}.  In this paper we use the original definition given in\cite{AS}.  That is, the matrix $u$ is \emph{degenerate} if there exist elements $\sigm, \sign \in G$ and elements $a,b \in K^*$ so that $\langle \sigm,\sign\rangle$ is noncyclic and $u_{\barm,\barn}=\sigm(a)a^{-1}\sign(b)b^{-1}$.  The motivating consequence behind this definition is that if $u$ is degenerate, then $bz^{\barm}$ and $a^{-1}z^{\barn}$ commute in $\Delta$.

Recall that if $\Delta$ is an abelian crossed product then the {\it generic abelian crossed product} associated to $\Delta$, which we will denote by $\mathcal{A}_{\Delta}$, is the abelian crossed product 
\[\mathcal{A}_{\Delta} \cong (K(x_1,\ldots,x_r)/F(x_1,\ldots,x_r),G,z,u,bx).\]
Here $x_1,\ldots,x_r$ are independent indeterminates and $bx=\{b_ix_i\}_{i=1}^r$.  Let $p$ be a prime.  A {\it $p$-algebra} is a central simple algebra over a field of characteristic $p$ with $p$-power index.  As the main result in \cite{AS}, generic abelian crossed products and non-degeneracy were used to establish the existence of non-cyclic $p$-algebras.  

\subsection{Related examples}
As mentioned above, the results in \cite[Theorem 2.3.1]{mck-galoissubfields} (in the case $\chara(F)=p$) and \cite[Cor. 3.6]{mounirh1} show that a generic abelian crossed product of $p$-power degree defined by a non-degenerate matrix is indecomposable.  In \cite[Section 3.3]{mck-galoissubfields} examples of such algebras with index $p^n$ and exponent $p$ for all $p \ne 2$ and all $n \geq 2$ are given.  An example is also given in the case $p=2$ and $n=3$.  In the $p\ne2$ example the abelian crossed product $\Delta$ is constructed by generically lowering the exponent of an abelian crossed product with exponent equal to index equal to $p^n$.  That is, let $\Delta'$ be an abelian crossed product defined by the group $G\cong (\zp)^n$, $n\geq 2$, with index and exponent $p^n$.  Set $Y=SB(\Delta'^{\otimes p})$ and let $\mathcal{F}(Y)$ be the function field of $Y$.  Then set $\Delta=\Delta'\otimes \mathcal{F}(Y)$.  $\Delta$ has exponent $p$ and index $p^n$ by the index reduction theorem of \cite{Schofield}.  $\Delta$ is shown to be defined by a non-degenerate matrix by studying the torsion in $\textrm{CH}^2(SB(\Delta))$ (\cite[Prop. 3.1.1]{mck-galoissubfields}).  Since $\Delta$ is defined by a non-degenerate matrix, $\mathcal{A}_{\Delta}$, the associated generic abelian crossed product, is indecomposable. 

Because of the method of construction of $\Delta$, by \cite[Corollary 5.4]{Karpenko2}, the abelian crossed product $\Delta$ is itself an indecomposable division algebra of exponent $p$ and index $p^n$.  In this paper we construct an abelian crossed product, $\Delta_{\F}$, which is decomposable of index $p^2$, exponent $p$ ($p\ne2$), and is defined by a non-degenerate matrix.  As mentioned above $\mathcal{A}_{\Delta_\F}$, the generic abelian crossed product associated to $\Delta_{\F}$, is therefore indecomposable.  The strategy is to make an abelian crossed product decomposable in a generic way and prove that the matrix defining the resulting decomposable abelian crossed product is non-degenerate.  The method in this paper does not use the results of \cite{Karpenko2}.  In particular, we avoid using the Chow group of Severi-Brauer varieties, hence providing a more elementary, though still technical, approach to decomposability of abelian crossed products.

\subsection{Outline of paper.}
In section \ref{section2} we construct $\Delta_{\F}$, a decomposable abelian crossed product division algebra of index $p^2$ and exponent $p$ (Lemma \ref{lemma4} and Corollary \ref{cor1}).  We study $\Delta_\F$ for the rest of the paper, with our goal being to show that it is defined by a non-degenerate matrix.  The difficulty in proving non-degeneracy of the  matrix defining $\Delta_{\F}$ lies in the fact that the lattice $M_{\omega}$ used in the definition of $\F$ is not $H^1$-trivial.  In section \ref{section3} we alleviate this problem by constructing an $H^1$-trivialization, $M$, of the lattice $M_\omega$ and analyzing its structure as a module over a group ring.  In section \ref{section4} we study the form of elements in $M$ which could possibly make the matrix degenerate.  Moreover, it is noted that it suffices to prove the matrix is non-degenerate in the lattice $M$.  
Finally in section \ref{section5} we prove the main theorem, Theorem \ref{p1}, which states that the matrix defining the abelian crossed product $\Delta_\F$ is non-degenerate.  

{\bf Acknowledgments} The author would like to thank David Saltman and Adrian Wadsworth for help with this project, especially for their help with the homological formulation of the degeneracy condition given in section \ref{section5}.

\section{The example}
\label{section2}
The goal of this section is to construct, in a generic way, a decomposable abelian crossed product of index $p^2$ and exponent $p$.  This is done using fields generated by group lattices as in \cite{Saltman-invariant}, \cite[section 12]{LN} and \cite[section 3.2]{mck-galoissubfields}.  We will recall the relevant objects as we need them in this section.

For any finite abelian group $H$ of rank $r$, generated by $\{\sigma_i\}_{i=1}^r$, let $I[H]$ be the augmentation ideal.  That is, $I[H]$ is the kernel $0 \to I[H] \to \Z[H] \stackrel{\epsilon}{\to} \Z \to 0$ where $\epsilon(\sigma)=1$ for all $\sigma \in H$.  Define $A_2(H)$ to be the $H$-lattice which is the kernel
\begin{equation}
0\longrightarrow A_2(H)\longrightarrow\bigoplus_{i=1}^r\Z[H] d_i\longrightarrow I[H] \longrightarrow0\\
\label{eqn2}\end{equation}
where $d_i$ is mapped to $\sigma_i-1$ for all $1\leq i\leq r$.  We note here that $A_2(H)$ is also known as the {\it relation module} for a minimal free presentation of $H$ (see e.g. \cite[Prop. 2.3]{Gruenberg}).  Let $[c_H]\in H^2(H,A_2(H))$ be the class of the so called ``canonical" 2-cocycle.  That is, if $\delta$ is the co-boundary map we set $[c_1]=\delta([1])\in H^1(H,I[H])$ for $[1]\in H^0(H,\Z)\cong \Z$.  Then, $[c_H]=\delta([c_1])$.  The module $A_2(H)$ and cocycle $[c_H]$ are also used in \cite[chapter 12]{LN} and \cite[section 3.2]{mck-galoissubfields}. 

For the rest of this paper let $p$ be a prime, $p \ne 2$, and for $i=1,2,3,4$, let $G_i=\langle \sigma_i \rangle \cong C_p$, the multiplicative cyclic group of order $p$.  Set $G=G_1\times G_2\times G_3\times G_4$ and fix the following notation:   for any numbers $i,j,k$ between 1 and 4, let $G_{ij}=G_i\times G_j$, and $G_{ijk}=G_i\times G_j \times G_k$,  each considered as a subgroup of $G$ .  For any subgroup $H \leq G$, define 
\begin{equation}
H_{12}=HG_{34}/G_{34},\hspace{.25in} H_3=HG_{124}/G_{124}\hspace{.25in} \textrm{and} \hspace{.25in} H_4=HG_{123}/G_{123}\label{eqn16}\end{equation}
so that $H_{12}\leq G_{12}$, $H_3\leq G_3$ and $H_4 \leq G_{4}$.  Furthermore, to ease the notation slightly, set $[c_{12}]=[c_{G_{12}}]$, $[c_3]=[c_{G_3}]$ and $[c_4]=[c_{G_4}]$.  These are the cocycle classes we use to build our algebra.

Let $F$ be a field.  Set $L_{12}=F(\atg)=q(F[\atg])$, the field of fractions of the commutative group ring $F[\atg]$ which is a domain.  The trivial $G_{12}$-action on $F$ and the natural $G_{12}$-action on $\atg$ extend to a $G_{12}$-action on $L_{12}$.  Since the $G_{12}$-action on $\atg$ is faithful, $L_{12}/L_{12}^{G_{12}}$ is a $G_{12}$-Galois extension of fields.  

For any group $H$, $H$-lattice $\Lambda$ and field $K$, let $e:\Lambda \to K[\Lambda]$ denote the canonical injection taking the additive group $\Lambda$ to the multiplicative subgroup of $K[\Lambda]$ consisting of monomials with coefficient 1.  By \cite[12.4(a)]{LN}, the associated map on cohomology $H^2(H,\Lambda) \to H^2(H,K(\Lambda)^*)$ is an injection and as such we will not distinguish between cocycle classes in $H^2(H,\Lambda)$ and their image in $H^2(H,K(\Lambda)^*)$.   

For $i=1,2,3,4$ set $\norm_i=1+\sigma_i+\ldots+\sigma_i^{p-1}$ and for $i,j \in \{1,2\}$, set 
\begin{eqnarray*}
b_i &=& \norm_i{d_i}\in\atg \\
u_{ij} &=& (\sigma_i-1)d_j-(\sigma_j-1)d_i\in\atg.\end{eqnarray*}
Define $e(u)=(e(u_{ij}))\in M_2(L_{12}^*)$ and $e(b)=(e(b_1),e(b_2))\in (L_{12}^*)^{2}$.  The matrix $e(u)$ and the vector $e(b)$ satisfy the conditions in \cite[Theorem 1.3]{AS}, and therefore 
\begin{equation}\Delta_{12}=(L_{12}/L_{12}^{G_{12}},z_{\sigma},e(u),e(b))\label{eqn3}\end{equation}
is an abelian crossed product.  Furthermore, as noted in \cite[Lemma 3.2.3]{mck-galoissubfields}, there is an isomorphism $\Delta_{12}\cong(L_{12}/L_{12}^{G_{12}},G_{12},c_{12})$.

Now we use the groups $G_3$ and $G_4$ to construct two generic degree $p$ cyclic algebras, $\Delta_3$ and $\Delta_4$.  For $j=3,4$, there is an isomorphism $A_2(G_j) \cong \norm_j \Z$.  Set $K_j=A_2(G_j) \oplus \Z[G_{j}]$.  Set $L_j=F(K_j)$ where again we take $G_j$ to act trivially on $F$.  Since $G_j$ acts faithfully on $K_j$, $L_j/L_j^{G_j}$ is a cyclic Galois extension of degree $p$.  For $j=3,4$ set,
\begin{equation}
\Delta_j=\left(L_j/L_j^{G_j}, G_j, c_j\right)
\end{equation}
where $c_j$ is a 2-cocycle in the canonical class $[c_j]\in H^2(G_j,L_j^*)$.  One can show that as a cyclic algebra $\Delta_j=(F(\Z[G_j])(x_j)/F(\Z[G_j])^{G_j}(x_j), x_j)$ where $x_j$ stands for the element $N_j \in K_j$ (see the proof of Lemma \ref{lemma5}).

Set $Q=A(G_{12})\oplus K_3\oplus K_4$.  $Q$ is a $\zg$-lattice with $G$-action given as follows.  $G_{12}$ acts in the natural way on $A(G_{12})$ and trivially on $K_3$ and $K_4$.  $G_3$ (resp. $G_4$) acts trivially on $A(G_{12})$ and $K_4$ (resp. $A(G_{12})$ and $K_3$) and acts in the natural way on $K_3$ (resp.  $K_4$).  Set $L=F(Q)$.  The action of $G$ on $Q$ extends to $L$ and since $G$ acts faithfully on $Q$, $L/L^G$ is a $G$-Galois extension.  Using the three inclusions $L_{12}^{G_{12}}, L_3^{G_{3}}, L_4^{G_4} \subset L^G$ we have the following three field diagrams:
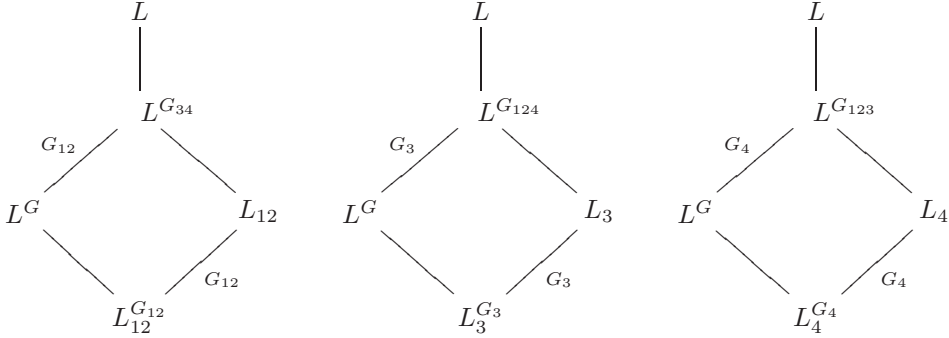
\[\xymatrix@C=10pt{
				&L			&				\\
				&{\hspace{.3in}}L^{G_{34}}\ar@{-}[u]			&				\\
L^G\ar@{-}[ur]^{G_{12}}	&					&L_{12}\ar@{-}[ul]	\\
				&L_{12}^{G_{12}}\ar@{-}[ur]_{G_{12}}\ar@{-}[ul]&\\
}\hspace{.25in}
\xymatrix@C=10pt{
				&L			&				\\
				&{\hspace{.3in}}L^{G_{124}}\ar@{-}[u]			&				\\
L^G\ar@{-}[ur]^{G_{3}}	&					&L_3\ar@{-}[ul]	\\
				&L_3^{G_{3}}\ar@{-}[ur]_{G_{3}}\ar@{-}[ul]&\\
				}\hspace{.25in}
\xymatrix@C=10pt{
				&L			&				\\
				&{\hspace{.3in}}L^{G_{123}}\ar@{-}[u]			&				\\
L^G\ar@{-}[ur]^{G_{4}}	&					&L_4\ar@{-}[ul]	\\
				&L_4^{G_{4}}\ar@{-}[ur]_{G_{4}}\ar@{-}[ul]&\\
				}
\]

We finally define $A$ to be the central simple $L^G$-algebra
\begin{equation}
A = (\Delta_{12} \otimes_{\lot}L^G)\otimes_{L^G} \left((\Delta_3\otimes_{L_3^{G_3}}L^G )\otimes_{L^G}(\Delta_4 \otimes_{L_4^{G_4}}L^G )\right)^{\circ},
\end{equation}
where ${}^\circ$ denotes the opposite algebra.  Set $[\omega] =[\ct]-[\ctn{3}]-[\ctn{4}]\in H^2(G,L^*)$, where $[\ct]\in H^2(G_{12},L_{12}^*)$, $[\ctn{3}]\in H^2(G_3,L_3^*)$ and $[\ctn{4}]\in H^2(G_4,L_4^*)$ are extended to $G$ by inflation.

\begin{lemma} $A \cong (L/L^G, \omega)$, where $\omega$ is a 2-cocycle in the class of $[\omega] \in H^2(G,L^*)$. 
\label{l1}\end{lemma}

\proof By \cite[(29.13) \& (29.16)]{Reiner} $\Delta_{12}\otimes L^G$ is similar to $(L/L^G,c_{12}')$, where $c_{12}'$ is a 2-cocycle in the image of the class $[c_{12}]$ under $\inf: H^2(G_{12},L_{12}^*) \to H^2(G,L^*)$.
Similarly, $\Delta_j \sim (L/L^G,c_j')$, where $c_j'$ is a 2-cocycle in the image of the class $[c_j]$ under
$\inf: H^2(G_{j},L_{j}^*) \to H^2(G,L^*)$.  Consequently, $A \sim (L/L^G,\omega)$.  Since the degree of $A$ equals the degree of $(L/L^G,\omega)$, this is an isomorphism.\endproof

Let $\F$ be the function field of $\textrm{SB}(A)$, the Severi-Brauer variety of $A$.  Since $A$ is a $G$-crossed product there is an explicit description of $\F$ as follows.  By \cite[Theorem 0.5]{Saltman-invariant}, $\F\cong F(Q)_{\omega}(I[G])^G=L_{\omega}(\ig)^G$.  Equivalently, $\F$ can be described by the following construction of the $G$-lattice $\mo$.  As an abelian group there is an isomorphism $\mo \cong Q \oplus \ig$ and the $G$-action on $\mo$ is defined by 
\begin{align}
g(x,0)&=(g\cdot x,0) \textrm{ for } x \in Q \textrm{ and }\nonumber\\
g(0,g'-1)&=(\omega(g,g'),g(g'-1)) \textrm{ for } g'-1 \in \ig.
\label{eqnlast}
\end{align}  
The field $\F$ satisfies the isomorphism $\F \cong F(\mo)^G$.

Since $\F$ is a splitting field for $A$ and the dimensions on both sides of \eqnref{eqn8} are equal, we have the isomorphism,
\begin{equation}
\Delta_{12} \otimes_{\lot} \F \cong (\Delta_3\otimes_{L_3^{G_3}} \F)\otimes_{\F} (\Delta_4 \otimes_{L_4^{G_4}} \F).\label{eqn8}
\end{equation}
In particular, $\Delta_{\F}:=\Delta_{12} \otimes_{\lot} \F$ is a decomposable abelian crossed product.  
\begin{lemma} $\Delta_{\F}\cong (F(\mo)^{G_{34}}/F(\mo)^G,G_{12},\ct)$ is a decomposable abelian crossed product. \label{lemma4}
\end{lemma}
\proof  
Since we have already noticed that $\Delta_\F$ is decomposable, we need only show the isomorphism.  $\Delta_{\F}$ is similar in the Brauer group to $(L_{12}\F/\F,H,\ct')$ where 
$\ct'$ is the restriction of $\ct$ to $H=\gal(L_{12}\F/\F)$.  
However, $L_{12} \cap F(\mo)^{G} =\lot$ and therefore $\gal(L_{12}\F/\F)=G_{12}$.  
Thus the given similarity is also an isomorphism.  Therefore, $L_{12}\F$ has degree $p^2$ over $\F=F(\mo)^{G}$ and is a maximal subfield of $\Delta_{\F}$.  Since both $L_{12}$ and $\F=F(\mo)^G$ are contained in $F(\mo)^{G_{34}}$, a degree $p^2$ field extension over $F(\mo)^G$, the composite must satisfy $L_{12}\F \cong F(\mo)^{G_{34}}$.  \endproof

By Lemma \ref{lemma4} and equation \eqnref{eqn3}, there is an isomorphism 
\begin{equation}\Delta_{\F} \cong (F(\mo)^{G_{34}}/F(\mo)^G,G_{12},z,e(u),e(b)).\label{eqn4}\end{equation} 
\begin{corollary}
$\Delta_{\F}\cong (F(M_\omega)^{G_{124}}/F(M_\omega)^G,x_3)\otimes (F(M_\omega)^{G_{123}}/F(M_\omega)^G,x_4)$ and as a consequence, there is a maximal abelian subfield $N_1\subset \Delta_{\F}$ with $N_1 \cong F(M_\omega)^{G_{12}}$ and $\Delta_{\F} \in \dec(N_1/F(M_\omega)^G)$.
\label{cor2}
\end{corollary}
\proof A similar argument as in the proof of Lemma \ref{lemma4} shows that $\Delta_3\otimes \mathcal F_A \cong (F(M_\omega)^{G_{124}}/F(M_\omega)^G,x_3)$ and $\Delta_4 \cong (F(M_\omega)^{G_{123}}/F(M_\omega)^G,x_4)$.  Since $F(M_\omega)^{G_{124}}\otimes F(M_\omega)^{G_{123}} \cong F(M_\omega)^{G_{12}}$ we see $\Delta_{\F}$ is split by $F(M_\omega)^{G_{12}}$.  Since $[F(M_\omega)^{G_{12}}:F(M_\omega)^G]=p^2=\deg(\Delta_{\F})$, there exists a maximal abelian subfield $N_1\subset \Delta_{\F}$ such that $N_1\cong F(M_\omega)^{G_{12}}$ and $\Delta_{\F} \in \dec(N_1/F(M_\omega)^G)$.
\endproof

We can now state the main theorem of the paper.  The proof is given in section \ref{section5}.
\begin{theorem} Let $F$ be a field with a $G$-action so that $F^*$ is an $H^1$-trivial $G$-module.  Then, \[\Delta_{\F}\cong \left(F(\mo)^{G_{34}}/F(\mo)^G,G_{12},z,e(u),e(b)\right)\] is a decomposable abelian crossed product division algebra defined by a non-degenerate matrix of index $p^2$ and exponent $p$, $p \ne 2$.\label{t1}
\end{theorem}

\begin{remark}
By Lemma \ref{lemma4} $\Delta_\F$ is decomposable, hence to prove this theorem there are two things left to show.  First we need to show that $\Delta_{\F}$ has index $p^2$ and therefore is a division algebra and second that $e(u)$ is non-degenerate in $F(\mo)^{G_{34}}$.  The difficulty will lie in showing that $e(u)$ is non-degenerate.  This is difficult because the way things stand now, there is no ``easy to understand" $G$-action on $F(\mo)$.  Recall that a $G$-lattice $\Lambda$ is said to be {\it $H^1$-trivial} if $H^1(H,\Lambda)=0$ for all subgroups $H\leq G$.  By \cite[Theorem 12.4(c)]{LN}, if $\Lambda$ is $H^1$-trivial, and $K$ is a field with $G$-action so that $K^*$ is also an $H^1$-trivial $G$-module, then $K(\Lambda)^*\cong_G K^*\oplus\Lambda \oplus P$, where $P$ is a permutation module and the isomorphism is a $G$-module isomorphism.    In particular, if the lattice $\Lambda$ is $H^1$-trivial, the $G$-action is ``easy to understand".  In section \ref{section3} we show that $\mo$ is not an $H^1$-trivial $G$-module.  As a consequence in section \ref{section3} we also construct an $H^1$-trivialization of $\mo$.  In Theorem \ref{p1} we show that $e(u)$ is non-degenerate in the field extension generated by this larger lattice. It is in the proof of the non-degeneracy of $e(u)$ that we use $p\ne 2$ (see the proof of Lemma \ref{lemma3}).
\label{remark5}
\end{remark}
\begin{remark}
In Theorem \ref{t1}, we can always choose the field $F$ so that it is a trivial $G$-module.  For example, we could choose $F$ to be a field of characteristic $p$ and then, since $G$ is a $p$-group, $F^*$ with trivial $G$-action is a trivial $H^1$-module.
\end{remark}
\begin{remark}
Theorem \ref{t1} shows that the converse to \cite[Prop. 3.1.1]{mck-galoissubfields} is false.  That is, \cite[Prop. 3.1.1]{mck-galoissubfields} says that if $\Delta$ is an abelian crossed product with exponent $p$ defined by a non-cyclic group $G$ and $u$ is degenerate ($p \ne 2$), then $\textrm{CH}^2(SB(\Delta))$ is torsion free.  In our case the algebra in Theorem \ref{t1} has exponent $p$.  Moreover $\textrm{CH}^2(SB(\Delta_{\mathcal{F}_A}))$ is torsion free since $\Delta_{\mathcal{F}_A}$ is decomposable (see \cite{Karpenko4} and \cite{Karpenko2}).  However, the matrix defining $\Delta_{\mathcal{F}_A}$ is non-degenerate.  The converse to \cite[Prop. 3.1.1]{mck-galoissubfields}, therefore, is false.
\label{remark3}
\end{remark}
\begin{remark}
Since $e(u)$ is a non-degenerate matrix, $\Delta_{\F} \notin \dec(F(M_\omega)^{G_{34}}/F(M_\omega)^G)$ (see e.g., \cite[Prop. 7.13]{LN}).  We also record this result in Corollary \ref{cor3}.
\label{remark4}
\end{remark}

\subsection{Index Calculation}
By the Schofield-Van den Bergh index reduction formula \cite{Schofield}, $\ind(\Delta_{12}\otimes \F)$ is the minimum of the index of $\Delta_{12}\otimes A^i$ as $i$ varies and the tensor product is taken over $L^G$.  To get a lower bound on this index, we can compute the exponent of this algebra which is precisely the exponent of the cocycle class $(i+1)[c_{12}]-i[c_3]-i[c_4] \in H^2(G,Q)$.  It is easy to see that
\[\exp(\Delta_{12}\otimes A^i)=\exp((i+1)[c_{12}]-i[c_3]-i[c_4])=\left\{\begin{array}{ll} p&p\mid i+1\\p^2&p\nmid i+1\end{array}\right..\]
Therefore, to show that $\ind(\Delta_{12}\otimes \F)=p^2$ one need only show when $p\mid i+1$, $\ind(\Delta_{12}\otimes A^i)=p^n$ for some $n\geq 2$.  This is accomplished by the following lemma.

\begin{lemma}Set 
\[B_m=(\Delta_{12}\otimes_{L_{12}^{G_{12}}}L^G)^{\otimes mp}\otimes_{L^G} (\Delta_3\otimes_{L_3^{G_3}}L^G)\otimes_{L^G}(\Delta_4\otimes_{L_4^{G_4}}L^G).\]
Then for all $m \in \Z^+$, $\ind(B_m)=p^n$ for some $n\geq 2$.\label{lemma5}
\end{lemma}

\proof  The idea of the proof is to show that the restriction of $B_m$ to $L^{G_{34}}$, a maximal subfield of $\Delta_{12}\otimes L^G$, has index $p^2$.  Set $E=F(\atg\oplus\Z[G_3]\oplus\Z[G_4])$.  $E$ is naturally a subfield of $L=F(Q)$ and in fact $L=E(x_3,x_4)$ where $x_3$ and $x_4$ are independent indeterminates corresponding to the factor of $\norm_3\Z$ and $\norm_4\Z$ in the lattice $Q$.  Moreover, $x_3$ and $x_4$ have trivial action by $G$ and therefore, $L^G=E^G(x_3,x_4)$.  The algebras $\Delta_3$ and $\Delta_4$ are the cyclic algebras $(F(\Z[G_j])(x_j)/F(\Z[G_j])^{G_{j}}(x_j),x_j)$ for $j=3,4$ respectively.  We have, 
\begin{align*}\Delta_3\otimes_{L_3^{G_3}}L^G \cong (E^{G_{124}}(x_3,x_4)/E^G(x_3,x_4),x_3),\\
\Delta_4\otimes_{L_4^{G_4}}L^G \cong (E^{G_{123}}(x_3,x_4)/E^G(x_3,x_4),x_4)
\end{align*} and 
\[\Delta_{12}\otimes_{L_{12}^{G_{12}}}L^G \cong (E^{G_{34}}(x_3,x_4)/E^G(x_3,x_4),c_{12}).\]
Therefore, 
\begin{multline}
B_m\otimes_{L^G} E^{G_{34}}(x_3,x_4)\sim\\
 (E^{G_4}(x_3,x_4)/E^{G_{34}}(x_3,x_4),x_3)\otimes (E^{G_3}(x_3,x_4)/E^{G_{34}}(x_3,x_4),x_4).
\label{eqn5} \end{multline}
The algebra in \eqnref{eqn5} is isomorphic to the central localization of the iterated twisted polynomial ring $E^{G_{34}}[t_3,t_4;\sigma_3,\sigma_4]$ which is a domain (see e.g., \cite[pg. 9]{LN}).  Therefore, $(E^{G_4}(x_3,x_4)/E^{G_{34}}(x_3,x_4),x_3)\otimes (E^{G_3}(x_3,x_4)/E^{G_{34}}(x_3,x_4),x_4)$ has index $p^2$.  This proves the lemma.\endproof

\begin{corollary}
$\Delta_{\F}$ has exponent $p$ and index $p^2$.
\label{cor1}\end{corollary}
\proof By the exponent calculation before Lemma \ref{lemma5} we only need to show that $p^2 \mid \ind(\Delta_{12}\otimes A^i)$ when $p|i+1$.  Assume $p\mid i+1$ and set $pm=i+1$.  Since $\Delta_3$ and $\Delta_4$ have index $p$, $\Delta_{12}\otimes A^i \sim B$ and hence the corollary follows from Lemma \ref{lemma5} and the fact that $\ind(\Delta_{\F})\mid \ind(\Delta_{12})=p^2$.\endproof

\section{An $H^1$-trivialization of $\mo$}
\label{section3}
Let $\mo$ be the $G$-lattice constructed in section \ref{section2}, equation \ref{eqnlast}.  In this section we construct a specific $H^1$-trivialization of $\mo$.  In particular, we construct an extension of $G$-lattices $0\to \mo \to M \to P \to 0$ such that $P$ is a permutation lattice and $M$ is $H^1$-trivial.  As noted in Remark \ref{remark5}, we need this $H^1$-trivialization because it gives us an easier to understand $G$-action on the field $F(M)$.  In particular, using this specific $H^1$-trivialization, $M$, in Theorem \ref{p1} we are able to show that $e(u)$ is non-degenerate in $F(M)^{G_{34}}$.  This implies that $e(u)$ is non-degenerate in $F(\mo)^{G_{34}}$ because of the natural inclusion $F(\mo) \subset F(M)$.  In fact, by \cite[Theorem 12.9]{LN}, $F(M)^{G_{34}}$ is a rational extension of $F(\mo)^{G_{34}}$ since $M$ is an extension of $\mo$ by a permutation lattice, however we do not need or use this fact.  

\begin{remark}Using the notation and terminology of \cite[Sections 2.5-2.7]{Lorenz}, one could define an $H^1$-trivialization of a $G$-lattice $\Lambda$ to be any $G$-lattice $\widetilde \Lambda$ such that $\widetilde \Lambda$ is $H^1$-trivial (or {\it coflasque} in \cite[2.5]{Lorenz}) and such that $\Lambda \sim_{\mathrm{fl}}\widetilde \Lambda$, that is, $\Lambda$ and $\widetilde \Lambda$ are {\it flasque equivalent}.  In particular, using \cite[2.7.1(c)]{Lorenz} we see that $\widetilde \Lambda$ is an $H^1$-trivialization of $\Lambda$ if and only if $\widetilde \Lambda$ is $H^1$-trivial and there is a short exact sequence of $G$-lattices $0\to \Lambda \to \widetilde \Lambda \oplus P \to Q \to 0$ with $P$ and $Q$ permutation lattices. By \cite[12.5]{LN} $H^1$-trivializations exist.   It was pointed out to the author that from \cite[2.7.1(c)]{Lorenz} it follows that any two $H^1$-trivializations are {\it stably permutation equivalent}.  That is, if $\widetilde \Lambda$ and $\Lambda'$ are both $H^1$-trivializations of $\Lambda$, then there exist permutation lattices $P_1$ and $P_2$ such that $\widetilde \Lambda \oplus P_1\cong \Lambda'\oplus P_2$.  In particular, an $H^1$-trivialization is unique up to stable permutation equivalence.  In this section we take the time to construct a specific $H^1$-trivialization because we need to know in Lemma \ref{lemma2} the form of the specific 1-cocycles that we split.  The specific form of the 1-cocycles allows us to prove Proposition \ref{l5}, giving us a valuable tool for proving that $e(u)$ is non-degenerate in $F(M)^{G_{34}}$ in Theorem \ref{p1}.
\end{remark}

$\mo$ is not itself an $H^1$-trivial $G$-lattice and in the first lemma we calculate the groups $H^1(H,\mo)$ for all subgroups $H \leq G$.  

\begin{lemma} Let $H$ be a subgroup of $G$.  Then $H^1(H,\mo)\cong \Z/n\Z$, where 
\[n=\frac{|H|}{|\res^G_H(\omega)|}\]
where $|\res^G_H(\omega)|$ is the order of the class of $\res^G_H(\omega)$ in $H^2(H,Q)$.
\label{l4}
\end{lemma}
\proof  By \cite[Lemma 3.2.4]{mck-galoissubfields}, $\atg$ is $H^1$-trivial.  Since $A_2(G_j) \cong \Z \norm_j$, a trivial permutation lattice, $K_3$ and $K_4$ are each permutation lattices and therefore they are $H^1$-trivial.  Together this implies that $Q=\amm$ is $H^1$-trivial.  From the short exact sequence $0 \to Q \to \mo \to \ig \to 0$ we get the long exact sequence of cohomology,
\[ \ldots \to 0 \to H^1(H,\mo) \to H^1(H,\ig) \to H^2(H,Q) \to \ldots\]
It is easy to see $H^1(H,\ig) \cong \Z/|H|\Z$ and this group is generated by the class of the 1-cocycle $d_H:H \to \ig$ given by $d_H(h)=h-1$ for all $h \in H$.  By the definition of $\mo$, $[d_H] \mapsto [\res^G_H(\omega)] \in H^2(H,Q)$.  Therefore, $H^1(H,\mo)$ is cyclic of order $n$ where $n=|H|/|\res^G_H(\omega)|$.  \endproof

\begin{lemma}Let $H$ be a subgroup of $G$.  Then,
\[|\res^G_H(\omega)|=\max\left\{|H_{12}|, |H_3|, |H_4|\right\},\]
where $H_{12}$, $H_3$ and $H_4$ are the quotient subgroups of $G_{12}$, $G_3$ and $G_4$ (respectively) defined in \eqnref{eqn16}.\label{lemma6}
\end{lemma}
\proof  Let $E\leq G_{12}$ be a subgroup.  Then it is easy to show $|\res^{G_{12}}_E([\ct])|=|E|$.  This follows from looking at a segment of the long exact sequence of cohomology 
\[\ldots \to H^1(E,\Z[G_{12}]) \to H^1(E,I[G_{12}]) \to H^2(E,A(G_{12}))\to\ldots \]
Here $H^1(E,\Z[G_{12}])=0$ since $\Z[G_{12}]$ is a permutation module and by \cite[12.3]{LN} all permutation modules are $H^1$-trivial.  The class, $\res^{G_{12}}_E([\ct])$ is the image of the generator of $H^1(E,I[G_{12}])$, a cyclic group of order $|E|$.  

Let $H$ be a subgroup of $G$.  Since $H$ is a $p$-group and $[\omega]=[\ct]-[\ctn{3}]-[\ctn{4}]$, the order of $\res^G_H([\omega])$ is the maximum of the orders of $[\ct]$, $[\ctn{3}]$ and $[\ctn{4}]$, inflated to $G$ and then restricted to $H$. In particular, by the above paragraph, $|\res^G_H([\ct])|=|H_{12}|$.  For $j=3$ and $4$, we have $|\res^G_{G_j}([c_j])|=p$.  Since $G_j$ has no non-trivial subgroups, we see that $|\res^G_H([c_j])|=|H_j|$.  The lemma now follows directly.\endproof 

Combining Lemmas \ref{l4} and \ref{lemma6} we immediately get the following corollary.
\begin{corollary}
$H^1(H,\mo) \cong \Z/n\Z$ where 
\[n=\frac{|H|}{\max\left\{|H_{12}|,|H_3|,|H_4|\right\}}.\]
\label{c1}
\end{corollary}
By \cite[12.5]{LN} one can construct an $H^1$-trivial module containing $\mo$ by splitting all non-trivial cocycles with permutation modules.  We will proceed in a slightly more efficient manner by splitting the non-trivial cocycles from the following subset $\mathcal{H}$ of subgroups of $G$.

Let $\mathcal{H}=\{G,\{\langle\tau,\sigma_3,\sigma_4\rangle|\tau \in G_{12}\}\}$, a set of subgroups of $G$.  For all $H \in \mathcal{H}$, let $f_H$ be a 1-cocycle whose class $[f_H]\in H^1(H,\mo)$ is mapped to $|\res^G_H(\omega)|[d_H]\in H^1(H,\ig)$.  In other words, by the proof of Lemma \ref{l4}, $[f_H]$ generates $H^1(H,\mo)$.  Define the $G$ module $M$ by the abelian group isomorphism 
\begin{equation}M \cong \mo \oplus P\label{e1}
\end{equation}
where $P$ is the permutation lattice $\bigoplus_{H \in \mathcal{H}} \Z[G/H]$.  For each $H \in \mathcal{H}$ fix a set of coset representatives $\{g_i\}$ for $G/H$ and let $\Z[G/H]$ be generated over $\Z$ by the symbols $u_{g_iH}$.  Define the $G$-action on $M$ by $g(x,0)=(g\cdot x,0)$ for $x \in \mo$ and $g(0,u_{g_iH})=(g_j f_H(h),u_{g_jH})$, where $gg_i=g_jh$ with the $g_i,g_j$ elements of the fixed coset representatives of $G/H$ and $h \in H$.  
\begin{lemma}$M$ is $H^1$-trivial.\label{lemma7}
\end{lemma}
\proof  From the short exact sequence $0 \to \mo \to M \to P \to 0$ and the fact that permutation modules are $H^1$-trivial, there is a surjection $H^1(K,\mo) \to H^1(K,M)$ for all subgroups $K$ of $G$.  Therefore, if we show that the generators of $H^1(K,\mo)$ are split in $M$, then we will have shown that $M$ is $H^1$-trivial.  By construction of the module $M$, for subgroups $H \in \mathcal{H}$, $H^1(H,M)=0$ because $(f_H(h),0)=(h-1)(0,u_H)$ for all $h \in H$.  Moreover, if the generator $[f_K]\in H^1(K,\mo)$ satisfies $[f_K]=\res^H_K([f_H])$ with the generator $[f_H] \in H^1(H,\mo)$ for some $H \in\mathcal{H}$, then $[f_K]$ is split in $M$ and therefore $H^1(K,M)=0$.  We will use these arguments in the cases below.

Let $K\leq G$ a subgroup of $G$ with $H^1(K,\mo) \cong \Z/n\Z$ with $n \ne 1$.  By Corollary \ref{c1} there are no $K\leq G$ with $n=|K|$.   We can also assume $K \ne G$ since $G \in \mathcal{H}$.  We address the three remaining cases separately.  \\
{\bf Case 1:}  $|K|=p^3$ and $n=p$.  In this case $|\res^G_K(\omega)|=p^2$ and therefore the kernel of $H^1(K,\ig) \to H^1(K,Q)$ is generated by $p^2[d_K]$.  Since restriction commutes with the long exact sequence of cohomology we have the following commutative diagram.
\[\xymatrix{
0\ar[r]&H^1(G,\mo)\ar[r]\ar[d]^{\res^G_K}&  H^1(G,\ig)\ar[r]\ar[d]^{\res^G_K}&  H^2(G,Q)\ar[d]^{\res^G_K}\\
0\ar[r]&H^1(K,\mo)\ar[r]&                       H^1(K,\ig)\ar[r]&                       H^2(K,Q)
}\]
Since $|[\omega]|=p^2$, $[f_G]\mapsto p^2[d_G]$.  Moreover since $\res^G_K(p^2[d_G])=p^2[d_K]$, the class $\res^G_K([f_G])\mapsto p^2[d_K]$.  Therefore, $\res^G_K([f_G])=[f_K]$ and since $[f_G]$ is split in $M$, so is $\res^G_K([f_G])$.  Therefore, $H^1(K,M)=0$.\\
{\bf Case 2:} $|K|=p^3$ and $n=p^2$.  In this case $|\res^G_K(\omega)|=p$.  We will show that $K\in \mathcal{H}$.  Since $n=p^2$, by Corollary \ref{c1}, $|K\cap G_{34}|= p^2$, and in particular, $G_{34} \leq K$.  Then $K/G_{34}$ is isomorphic to a cyclic $p$-subgroup of $G_{12}$, say generated by $\tau$.  Then, $K=\langle\tau,\sigma_3,\sigma_4\rangle \in \mathcal{H}$ and hence $H^1(K,M)=0$.\\
{\bf Case 3:} $|K|=p^2$ and $n=p$.  In this case $|\res^G_K(\omega)|=p$ and therefore $p[d_K]$ generates kernel of $H^1(K,\ig) \to H^2(K,Q)$.  We will show that $K$ is a subgroup of an $H \in \mathcal{H}-\{G\}$ and $\res^H_K([f_H])=[f_K]$.  Since $\max\{|K_{12}|,|K_3|,|K_4|\}=p$ and $|K_{12}|=|K|/|K\cap G_{34}|$, we see $|K\cap G_{34}|\geq p$.  Therefore there exist two elements, $\tau_1 \in G_{34}\cap K$ and $\tau_2\in K$ so that $K=\langle \tau_1,\tau_2\rangle$.  Write $\tau_2=\sigma_1^{m_1}\sigma_2^{m_2}\sigma_3^{m_3}\sigma_4^{m_4}$.  Then clearly, $K\leq H=\langle\sigma_1^{m_1}\sigma_2^{m_2},\sigma_3,\sigma_4\rangle \in \mathcal{H}-\{G\}$.  As in case 2 we have the following commutative diagram.
\[\xymatrix{
0 \to H^1(H,\mo)\ar[r]\ar[d]^{\res^H_K}&  H^1(H,\ig)\ar[r]\ar[d]^{\res^H_K}&  H^2(H,Q)\ar[d]^{\res^H_K}\\
0 \to H^1(K,\mo)\ar[r]&                       H^1(K,\ig)\ar[r]&                       H^2(K,Q)
}\]
By Lemma \ref{l4} $|\res^G_H(\omega)|=p$.  Therefore, $[f_H]\mapsto p[d_H]\in H^1(H,\ig)$ and since $\res^H_K(p[d_H])=p[d_K]$, the class $\res^H_K([f_H])\mapsto p[d_K]$.  Therefore, $\res^H_K([f_H])=[f_K]$.
\endproof

In section \ref{section4}, Proposition \ref{l5} we will need an explicit description of a cocycle in the class of $[f_H]\in H^1(H,\mo)$ for all $H \in \mathcal{H}$.  In fact Lemma \ref{lemma2} is the only result from Section \ref{section3} we use in the remainder of the paper.  For any $H \in \mathcal{H}$ let $v_H:H_{12} \to \atg$ be the 1-cochain defined by $v_H(\bar{h}_1)=\sum_{\bar{h} \in H_{12}}c_{12}(\bar{h}_1,\bar{h})$ where $H_{12}$ is considered as a subgroup of $G_{12}$.  Since $H_{12}$ is a finite group we have $|H_{12}|\,[\ct]=0$ in $H^2(H_{12},\atg)$.  In fact one can check that the 1-cochain $v_H$ satisfies $|H_{12}|\,c_{12}=\delta v_H$ where $\delta$ is the co-boundary map.  We will use the cochain $v_H$ in the proof of the following lemma.

\begin{lemma}  For every $H \in \mathcal{H}$ there is a 1-cochain $z_H:H \to K_3\oplus K_4$ so that the 1-cochain $f_H:H \to \mo$ defined by
\[f_H(h)=(z_H(h)-\inf^{H_{12}}_Hv_H(h),|\res^G_H(\omega)|(h-1))\]
is a 1-cocycle whose image in $H^1(H,\ig)$ is $\p\,[d_H]$.  
\label{lemma2}
\end{lemma}

\proof If the given cochain is a cocycle, then its image in $H^1(H,\ig)$ is clearly $|\res^G_H(\omega)|\, [d_H]$.  Hence we need only show that the given cochain is a cocycle for some cochain $z_H$.

Recall $[\omega]\in H^2(G,Q)$ is defined by inflation of $c_{12}$, $c_3$ and $c_4$ from the three subgroups $G_{12}$, $G_3$ and $G_4$, viewed as quotient groups of $G$.  Since inflation commutes with restriction (\cite[Prop. 1.5.5]{Neukirch}), 
\begin{equation}[\res^G_H(\omega)]=[\inf^{H_{12}}_H(c_{12})]+[\inf^{H_3}_H(c_3)]+[\inf^{H_4}_H(c_4)]\label{eqn6}\end{equation}
Note that for every $H \in \mathcal{H}$, $\max\{|H_{12}|,|H_3|,|H_4|\}=|H_{12}|$.  Therefore, by Lemma \ref{l4}, $|\res^G_H(\omega)|=|H_{12}|$.  Moreover, for every $H \in \mathcal{H}$, $p$ divides $|H_{12}|$.  Since $|[c_3]|=|[c_4]|=p$, there are 1-cochains $\widehat{c}_3:H_3 \to K_3$ and $\widehat{c}_4:H_4 \to K_4$ so that $\delta \widehat{c}_3=\p c_3$ and $\delta\widehat{c}_4=\p c_4$.  Using the fact that inflation commutes with the co-boundary homomorphism (\cite[Prop. 1.5.2]{Neukirch}) and the 1-cochain $v_H$ given above satisfies $|H_{12}|\,c_{12}=\delta v_H$, we can explicitly write $|\res^G_H(\omega)|\res^G_H(\omega)$ as
\[|\res^G_H(\omega)|\res^G_H(\omega)=\delta(\inf^{H_{12}}_Hv_H+\inf^{H_3}_H\widehat{c}_3+\inf^{H_4}_H\widehat{c}_4).\]
Set $z_H:H \to K_3\oplus K_4$ to be the 1-cochain $z_H=-\inf^{H_3}_H\widehat{c}_3-\inf^{H_4}_H\widehat{c}_4$.  Then,
\begin{eqnarray*}
\delta f_H (h_1,h_2)&=& (\delta z_H(h_1,h_2)-\delta \inf^{H_{12}}_Hv_H(h_1,h_2)+\p\res^G_H(\omega),0)\\
&=&(0,0).
\end{eqnarray*}
The cochain $f_H$ is a cocycle since its co-boundary is zero.\endproof

\section{Elements of $M^{G_{34}}$}  
\label{section4}
Let $M$ be the $H^1$-trivialization of the $G$-lattice $\mo$ defined in \eqnref{e1}.  Let $M^{G_{34}}$ indicate the elements in $M$ which are fixed by the subgroup $G_{34}\leq G$.  In this section we construct a $G_{12}$-module homomorphism $\pi':M^{G_{34}} \to \Z/p\Z$ which will be used to distinguish elements of $M^{G_{34}}$.  As a first step we recall some facts about $\atg$ since it naturally sits in $M^{G_{34}}$ as a $G_{12}$-submodule.
\begin{remark}
The next two lemmas can be gleaned from \cite[pg. 49]{LN}.  
\end{remark}

\begin{lemma}$\atg$ is generated over $\Z[G_{12}]$ by $u_{12}=(\sigma_2-1)d_1-(\sigma_1-1)d_2$, $b_1=\norm_{1}d_1$ and $b_2=\norm_{2}d_2$.
\label{lemma1}
\end{lemma}

\begin{lemma}  The relations in $\atg$ are generated over $\zgot$ by $(\sigma_i-1)b_i=0$ for $i=1,2$, $(\sigma_2-1)b_1=\norm_1u_{12}$ and $-(\sigma_1-1)b_2=\norm_2u_{12}$.   In particular, let $x,y,z\in \zgot$ such that $xu_{12}+yb_1+zb_2=0$.  Then, $x = z'N_2+y'N_1$, $y=-(\sigma_2-1)y'+(\sigma_1-1)y''$ and $z=(\sigma_1-1)z'+(\sigma_2-1)z''$ for some $y'$, $y''$, $z'$, $z'' \in \zgot$.
\label{lemma11}
\end{lemma}

By Lemma \ref{lemma11} $\atg\cong\left(\zgot u_{12}\oplus\zgot b_1\oplus \zgot b_2\right)/R$ where
\[R=\left\{xu_{12}+yb_1+zb_2\left| \begin{array}{l}x=(z'N_2+y'N_1)\\ y=-(\sigma_2-1)y'+(\sigma_1-1)y''\\z=(\sigma_1-1)z'+(\sigma_2-1)z''\\\textrm{for some }y',y'',z',z'' \in \zgot \end{array}\right.\right\}\]
Let $(\epsilon,\epsilon,\epsilon):\zgot u_{12}\oplus\zgot b_1\oplus \zgot b_2 \to \Z\oplus\Z\oplus\Z$ be the map induced from the augmentation map $\epsilon:\zgot \to \Z$.  Since $(\epsilon,\epsilon,\epsilon)(R)\subset p\Z\oplus 0 \oplus 0$, there is an induced map $\atg \to \Z/p\Z\oplus\Z\oplus \Z$.  Let $\pi:\atg \to \Z/p\Z$ be the composition of this map with projection onto the first coordinate.  $\pi$ is a $G_{12}$-module homomorphism where $G_{12}$ acts trivially on $\Z/p\Z$.  We will use the map $\pi$ to distinguish elements of $\atg$ from one another.  Note that if we define $\atg$ to be a $G$-module by assuming that $G_{34}$ acts trivially on it, then $\pi:\atg \to \Z/p\Z$ is also a $G$-module homomorphism.  

In the next series of lemmas we prove that there is a $G$-module homomorphism $\pi':M^{G_{34}} \to \Z/p\Z$, extending $\pi:\atg \to \Z/p\Z$.  In order to prove this, we first need to calculate the value of $\pi$ at some particular elements of $\atg$.  Recall for $\barm=(m_1,m_2),\barn=(n_1,n_2) \in \N^2$ we set $z^{\barm}=z_1^{m_1}z_2^{m_2} \in \Delta_{12}$ where $\Delta_{12}$ is the $G_{12}$-abelian crossed product given in \eqnref{eqn3}.  Set $u_{\barm,\barn}\in \atg$ to be the unique element such that $e(u_{\barm,\barn})=z^{\barm}z^{\barn}(z^{\barm})^{-1}(z^{\barn})^{-1}$.  

\begin{lemma}
For all $\barm=(m_1,m_2),\barn=(n_1,n_2)\in\N^2$, $\pi(u_{\barm,\barn})=m_1n_2-m_2n_1+p\Z$.
\label{l2}\end{lemma}
\proof Let $s,t \in \N$.  From the rules of multiplication in $\Delta_{12}$, it is easy to calculate
\begin{eqnarray}
e(u_{(s,0),(0,t)})&=&z_1^sz_2^t(z_1^s)^{-1}(z_2^t)^{-1}\nonumber\\
&=&\prod_{i=0}^{s-1}\prod_{j=0}^{t-1}\sigma_1^i\sigma_2^j(e(u_{12}))
\end{eqnarray}
One can also check that $z^{\barm}z^{\barn}=\sigma_1^{m_1}(e(u_{(0,m_2),(n_1,0)}))\sigma_1^{n_1}(e(u_{(m_1,0),(0,n_2)}))$.  Therefore, 
\begin{equation}u_{\barm,\barn}=-\sigma_1^{m_1}\cdot\left(\sum_{j=0}^{n_1-1}\sum_{i=0}^{m_2-1}\sigma_1^i\sigma_2^j\,u_{12}\right)+\sigma_1^{n_1}\left(\sum_{j=0}^{n_2-1}\sum_{i=0}^{m_1-1}\sigma_1^i\sigma_2^j\,u_{12}\right).\label{eqn13}\end{equation}
Applying $\pi$ to $u_{\barm,\barn}$ using formula \eqnref{eqn13} and the fact that $\pi(g\cdot x)=\pi(x)$ for all $g \in G_{12}$ and all $x \in \atg$ we immediately get $\pi(u_{\barm,\barn})=m_1n_2-m_2n_1+p\Z$.\endproof

\begin{lemma}
Let $H \in \mathcal{H}$ and $v_H:H_{12} \to \atg$ be the 1-cochain defined before Lemma \ref{lemma2}.  For $p \ne 2$, $\pi(v_H(\bar{h}))=0$ for all $\bar{h} \in H_{12}$.
\label{lemma3}
\end{lemma}
\proof Let $H \in \mathcal{H}$ and $h \in H$ with image $\bar{h} \in H_{12}$.  By definition, $v_H(\bar{h})=\sum_{\bar{h}'\in H_{12}}c_{12}(\bar{h},\bar{h}')$.  Let $\bar{h}=\sigm=\sigma_1^{m_1}\sigma_2^{m_2}\in H_{12}$ and let $\sign=\sigma_1^{n_1}\sigma_2^{n_2} \in H_{12}$.  By the definition of $c_{12}$, $e(c_{12}(z_g,z_h))=z_gz_h(z_{gh})^{-1}$ for all $g,h \in G_{12}$.  Let $m_i+n_i =q_ip+r_i$ with $0 \leq r_i <p$ for $i=1,2$.  Then,
\begin{eqnarray}
e(c_{12}(\sigm,\sign))&=&z^{\barm}z^{\barn}(z_1^{r_1}z_2^{r_2})^{-1}\nonumber\\
&=&\sigma_1^{m_1}(u_{(0,m_2),(n_1,0)})z_1^{m_1+n_1}z_2^{m_2+n_2}(z_1^{r_1}z_2^{r_2})^{-1}\nonumber\\
&=&\sigma_1^{m_1}(u_{(0,m_2),(n_1,0)})b_1^{q_1}\sigma_1^{r_1}(b_2^{q_2})\label{eqn15}
\end{eqnarray}
Note that if $\bar{h}=1$, then $c_{12}(\bar{h},\bar{h}')=0$ for all $\bar{h}' \in H_{12}$. We assume from now on that $\bar{h} \ne 1$.  We now consider the two cases $H=G$ and $H \ne G$ separately.  First, if $H \ne G$, then $H_{12}=\langle \bar{h}\rangle$ with $\bar{h}=\sigm$.  For each $0 \leq i \leq p-1$ set $im_1=q_ip+r_i$ with $0\leq r_i <p$.  Then, 
\begin{eqnarray*}
\pi(u(\bar{h}))
&=& \sum_{i=0}^{p-1}\pi(c_{12}(\sigm,(\sigm)^i))\nonumber\\
&=&\sum_{i=0}^{p-1}\pi(\sigma_1^{n_1}(u_{(0,m_2),(r_i,0)}))\hspace{.5in}\textrm{(by \eqnref{eqn15})}\nonumber\\
&=&-\frac{p(p-1)}{2}m_1m_2 + p\Z.
\end{eqnarray*}
The last line follows from Lemma \ref{l2} and the fact that $r_i \equiv im_i$ for all $i$.  For all $p \ne 2$, $-\frac{p(p-1)}{2}m_1m_2\equiv 0 \mod{p}$.  Hence $\pi(u(\bar{h}))\equiv 0$ in this case.  Similarly, if $H=G$, then $H_{12}=G_{12}$ and
\begin{eqnarray}
\pi(u(\bar{h}))
&=&\sum_{i=0}^{p-1}\sum_{j=0}^{p-1}\pi(c_{12}(\sigm,\sigma_1^i\sigma_2^j))\nonumber\\
&=&\sum_{i=0}^{p-1}\sum_{j=0}^{p-1}\pi(\sigma^{m_1}(u_{(0,m_2),(i,0)}))\hspace{.5in}\textrm{(by \eqnref{eqn15})}\nonumber\\
&=&\sum_{i=0}^{p-1}\sum_{j=0}^{p-1}-m_2\,i \hspace{.5in} \textrm{(by Lemma \ref{l2})}\nonumber\\
&=&-pm_2\frac{p(p-1)}{2}+p\Z=p\Z.\nonumber
\end{eqnarray}
Therefore $\pi(u(\bar{h}))\equiv 0$ in the case $H=G$.\endproof

Lemmas \ref{l2} and \ref{lemma3} will be used to prove that $\pi$ extends to $\pi':M^{G_{34}} \to \Z/p\Z$ in Proposition \ref{l5}.  In the next lemma we express elements of $M$ as 3-tuples $(x,y,z) \in Q\oplus \ig \oplus P\cong M$.  Recall $Q =\amm$ and $P = \bigoplus_{H \in \mathcal{H}}\Z[G/H]$.  In order to show the $G$-morphism $\pi:\atg \to \Z/p\Z$ extends to $M^{G_{34}}$, we will need to show that given $(x,y,z)\in M^{G_{34}}$ the $\ig$ component of the direct sum, $y$, has a particular form.  This is the content of the following lemma.

\begin{lemma}  Let $(x,y,z) \in M$ and let $\norm_{34}=\sum_{g \in G_{34}}g\in \zg$ be the $G_{34}$ norm.   If $(x,y,z) \in M^{G_{34}}$, then 
\[y=\norm_{34}\, y_1-py_2\]
where $y_1\in \zgot$ and $y_2 \in \zg$ with $\epsilon(y_2)=p\epsilon(y_1)$.
\label{l6}
\end{lemma}

\proof Let $g \in G_{34}$ and $u_{g_iH}\in \Z[G/H]$ for some $H \in \mathcal{H}$.  Let $(0,0,u_{g_iH}) \in M$.  Since $g\in H$ we have $gu_{g_iH}=u_{g_iH}$.  By the definition of the $G$-action on $M$ and Lemma \ref{lemma2}
\begin{equation}g(0,0,u_{g_iH})=(g_iz_H(g)+\p\omega(g_i,g),\p g_i(g-1),u_{g_iH}).\label{eqn12}\end{equation}

Here we have used the fact that $v_H(\bar{g})=0$ since $g \in G_{34}$.  Using the fact that $p|\p$ for all $H \in \mathcal{H}$, line \eqnref{eqn12} shows for each $z \in P$ and each $g \in G_{34}$ there exists $x_{z,g} \in Q$ and $y_z \in \zg$ so that $g(0,0,z)=(x_{z,g},y_zp(g-1),z)$.  Similarly, for each $y \in \ig$ and $g \in G_{34}$ there exists $x_{y,g}\in Q$ so that $g(0,y,0)=(x_{y,g},gy,0)$.  Finally we can conclude that for $g \in G_{34}$,
\[(g-1)(x,y,z)=(gx+x_{y,g}+x_{z,g}-x,(y_zp+y)(g-1),0).\]
Therefore, if $(x,y,z) \in M^{G_{34}}$, then $(g-1)(y+py_z)=0$ for all $g \in G_{34}$.  This implies that $y=\norm_{34}\cdot y_{1}-py_z$ for some $y_1 \in \Z[G_{12}]$.  Since $y \in \ig$, $\epsilon(y)=p^2\epsilon(y_1)-p\epsilon(y_z)=0$.  Setting $y_2=y_z$ completes the proof.  \endproof

\begin{proposition}  For $p \ne 2$ there exists a $G$-module homomorphism $\pi':M^{G_{34}} \to \Z/p\Z$ extending $\pi:\atg \to \Z/p\Z$.  
\label{l5}
\end{proposition}

\proof The $G_{12}$-module homomorphism $\pi:\atg \to \Z/p\Z$ easily extends to a $G$-module homomorphism on $Q=\atg \oplus K_3 \oplus K_4$ by setting $\pi(x_0,x_1,x_2)=\pi(x_0)$.  This follows because $Q$ is a direct sum as  $G$-modules and $G_{34}$ acts trivially on $\atg$.  Define $\pi':M^{G_{34}} \to \zp$ by
\[\pi'(x,y,z)=\pi(x).\]
Since $G_{34}$ acts trivially both on $\zp$ and $M^{G_{34}}$, to show that $\pi'$ is a $G$-module homomorphism we need only show that $\pi'$ is a $G_{12}$-module homomorphism.  Let $j \in\{1,2\}$ and let $(x,y,z) \in M^{G_{34}}$.  To prove our proposition we need to show 
\[\pi'(\sigma_j(0,y,0))=0\hspace{.25in}\textrm{and}\hspace{.25in}\pi'(\sigma_j(0,0,z))=0.\]
By Lemma \ref{l6}, $y=\norm_{34}\, y_1-py_2$ for some $y_1\in\zgot$ and $y_2 \in \zg$.  Note that for any $\alpha \in \ig$, $\pi'(0,p\alpha,0)=p\pi'(0,\alpha,0)=0$.  This fact is used many times in the calculations below.  Note also that $y_2-py_1 \in\ig$ by Lemma \ref{l6}.  Therefore,
\begin{eqnarray*}
\pi'(\sigma_j(0,y,0))&=&\pi'(\sigma_j(0,y_1(\norm_{34}-p^2)-p(y_2-py_1),0))\nonumber\\
&=&\pi'(\sigma_j(0,y_1(\norm_{34}-p^2),0))
\end{eqnarray*}
Set $y_1=\sum_{g \in G_{12}}\alpha_g g$ with $\alpha_g \in \Z$.  Then,
\begin{eqnarray}
\pi'(\sigma_j(0,y,0))&=&\sum_{g \in G_{12}}\alpha_g\, \pi'\left(\sigma_j(0,\sum_{h \in G_{34}}(gh-1)-p^2(g-1),0)\right)\nonumber\\
&=&\sum_{g \in G_{12}}\alpha_g \,\pi'\left(\sum_{h \in G_{34}}\omega(\sigma_j,gh),\sum_{h \in G_{34}}(gh-1),0\right)\nonumber\\
&=&\sum_{g \in G_{12}}\alpha_g\,\pi\left(\sum_{h \in G_{34}}\omega(\sigma_j,gh)\right)\nonumber\\
&=&\sum_{g \in G_{12}}\alpha_g\,\pi\left(p^2c_{12}(\sigma_j,g)\right)=0
\end{eqnarray}
To show $\pi'(\sigma_j(0,0,z))=0$ for all $j \in \{1,2\}$ and all $z \in P$, we show that $\pi'(\sigma_j(0,0,u_{g_iH}))=0$ for all $H \in \mathcal{H}$ and all coset representatives $g_i \in G/H$.  Set $\sigma_jg_i=g_kh$ where $g_k$ is another fixed coset representative of $G/H$ and $h \in H$.  Set $\p=p^{\delta_H}$ and note that for all $H \in \mathcal{H}$, $\delta_H>0$.  Then,
\begin{eqnarray*}
\pi'(\sigma_j(0,0,u_{g_iH}))&=&\pi'(g_kf_H(h),u_{g_kH})\nonumber\\
&=&\pi'(g_k(z_H(h)-v_H(\bar{h}))+p^{\delta_H}\omega(g_k,h),p^{\delta_H} g_k(h-1),u_{g_kH})\nonumber\\
&=&\pi(g_kz_H(h)-g_kv_H(\bar{h})+p^{\delta_H}\omega(g_k,h))\nonumber\\
&=&\pi(-v_H(\bar{h}))=0\hspace{.5in} 
\end{eqnarray*}
The last two equalities follow from the fact that $z_H(h) \in K_1\oplus K_2$ (Lemma \ref{lemma2}) and $\pi(v_H(\bar{h}))=0$ for all $\bar{h} \in H_{12}$ (Lemma \ref{lemma3}).
\endproof

\section{$e(u)$ is non-degenerate in $F(\mo)^{G_{34}}$}
\label{section5}
Before we prove the main theorem of the paper we give a homological formulation of the degeneracy condition.  Let $\Delta=(K/F,G,z,u,b)$ be an abelian crossed product with $G$ any finite abelian group.  Any commutator $[x,y]$ in $\Delta$ has reduced norm 1 hence any commutator that lands in $K$ has $K/F$-norm 1.  Therefore, using the $K$-basis of $\Delta$ given by $\{z^{\sigma}|\sigma \in G\}$ we get a map,
\begin{gather*}
\varphi:G\times G \to H^{-1}(G,K^*)\\
(\sigma,\tau) \mapsto [z^{\sigma},z^{\tau}]
\end{gather*}
In this map we are using the identification $H^{-1}(G,K^*) \cong \{N(K)=1\}/\ig K^*$  where $\{N(K)=1\}$ are the set of elements in $K^*$ with $K/F$-norm equal to 1 and $\ig K^*$ are the elements of the form $g(k)/k$ for $k \in K^*$ and $g \in G$.  Let $c(\sigma,\tau)=z^{\sigma}z^{\tau}(z^{\sigma\tau})^{-1}$ be the 2-cocycle associated to the abelian crossed product $\Delta$.  Using the properties of $H^{-1}(G,K^*)$ and the standard commutator identities for products we can show that $\varphi$ is bimultiplicative.
\begin{align*}
\varphi(\sigma\tau,\gamma)&=[z^{\sigma\tau},z^{\gamma}]=[c(\sigma,\tau)^{-1}z^{\sigma}z^{\tau},z^{\gamma}]\\
&=[z^{\sigma}z^{\tau},z^{\gamma}]\\
&=\sigma([z^{\tau},z^{\gamma}])[z^{\sigma},z^{\gamma}]\\
&=\varphi(\tau,\gamma)\varphi(\sigma,\gamma)
\end{align*}
(The identity $\varphi(\sigma,\tau\gamma)=\varphi(\sigma,\tau)\varphi(\sigma,\gamma)$ is done in an identical fashion).  Since $\varphi$ is also clearly symplectic, there is an induced map, which we also call $\varphi$, 
\[\varphi:G\wedge G \to H^{-1}(G,K^*)\]
\begin{lemma}  Let $\Delta$ be the abelian crossed product given above.  Then $u$ is degenerate if and only if there exists a rank 2 subgroup $H \leq G$ so that $\varphi|_H:H\wedge H \to H^{-1}(H,K^*)$ is the trivial map.
\label{lemma}
\end{lemma}
\proof Assume that $u$ is degenerate.  Then by definition there exits $\sigm,\sign \in G$ so that $H=\langle \sigm,\sign\rangle$ is a non-cyclic group and $u_{\barm,\barn}=\sigm(a)a^{-1}\sign(b)b^{-1}.$  In other words, $\varphi(\sigm,\sign)=u_{\barm,\barn}=0$ as an element of $H^{-1}(H,K^*)$.  Now by the bimultiplicativity of $\varphi$, $\varphi((\sigm)^s,(\sign)^t)=\varphi(\sigm,\sign)^{st}=0$.  Therefore, $\varphi|_H$ is the trivial map.  Conversely assume there is a rank 2 subgroup $H\leq G$ so that $\varphi|_H$ is the trivial map.  Set $H =\langle \sigm,\sign\rangle$.  Then $\varphi(\sigm,\sign) =0$, in other words, $u_{\barm,\barn}\in I[H] K^*$.  We are done since $I[H]$, the augmentation ideal of $\Z[H]$, is generated by $(\sigm-1)$ and $(\sign-1)$.
\endproof

In the case under investigation in this paper the group $G$ in the abelian crossed product is isomorphic to $\Z/p\Z\times \Z/p\Z$ and hence $G \wedge G \cong \Z/p\Z$.  Therefore to show that $u$ is non-degenerate in this case it suffices to show that $\varphi$ is non-trivial on a single commutator.  In particular, one need only show that $u_{12} \ne (\sigma_1)(a)a^{-1}\sigma_2(b)b^{-1}$ for all $a, b \in K^*$.

We can now prove the main theorem of the paper which states that the matrix defining the decomposable abelian crossed product $\Delta_\F$ is non-degenerate.  As always we take $p$ a prime, $p\ne 2$, $G=\langle \sigma_1\rangle\times \langle \sigma_2\rangle\times \langle \sigma_3\rangle\times\langle\sigma_4\rangle$ the elementary abelian group of order $p^4$, $G_{34}=\langle \sigma_3\rangle\times \langle \sigma_4\rangle$ and $G_{12}=\langle \sigma_1\rangle \times \langle \sigma_2\rangle$.

\begin{theorem}[Also listed as Theorem \ref{t1}] Let $F$ be a field with a $G$-action so that $F^*$ is an $H^1$-trivial $G$-module.  Then, \[\Delta_{\F}\cong \left(F(\mo)^{G_{34}}/F(\mo)^G,G_{12},z,e(u),e(b)\right)\] is a decomposable abelian crossed product division algebra defined by a non-degenerate matrix of index $p^2$ and exponent $p$, $p \ne 2$.
\label{p1}
\end{theorem}

\proof By Lemma \ref{lemma4} and Corollary \ref{cor1} $\Delta_\F$ is a decomposable abelian crossed product division algebra with isomorphism as stated in the theorem.  It is only left to show that $e(u)$ is a non-degenerate matrix.  As mentioned above the statement of the theorem, we need only show that there do not exist $a$, $b \in F(\mo)^{G_{34}}$ such that $e(u_{12})=\sigma_1(a)a^{-1}\sigma_2(b)b^{-1}$.  Since the $G$-lattice $\mo$ is a direct summand of $M$, there is an inclusion of fields $F(\mo)^{G_{34}} \subset F(M)^{G_{34}}$.  Therefore, it is enough to show that there do not exist $a$, $b \in F(M)^{G_{34}}$ such that $e(u_{12})=\sigma_1(a)a^{-1}\sigma_2(b)b^{-1}$.  By Lemma \ref{lemma7} $M$ is an $H^1$-trivial $G$-module.  Therefore, by \cite[12.4(c)]{LN},  we have a $G$-module isomorphism $F(M)^* \cong F^*\oplus M \oplus P'$.  Here $P'$ is a permutation $G$-module.  Under this isomorphism, $e(u_{12})\mapsto u_{12} \in \atg \subset M$.  Hence it is enough to show there do not exist $a,b \in M^{G_{34}}$ so that 
\begin{equation}
u_{12}=(\sigma_1-1)a+(\sigma_2-1)b.\label{eqn14}
\end{equation}
By Proposition \ref{l5} the map $\pi:\atg \to \zp$ extends to a $G$-module homomorphism $\pi':M^{G_{34}}\to \Z/p\Z$.  Taking $\pi'$ of both sides of \eqnref{eqn14} we get 
\[\pi'(u_{12})=1=\pi'((\sigma_1-1)a+(\sigma_2-1)b)=0.\]
This is a contradiction.  \endproof

\begin{corollary}
$\Delta_{\F} \notin \dec(F(M_\omega)^{G_{34}}/F(M_\omega)^G)$.
\label{cor3}
\end{corollary}
\proof By \cite[Prop. 7.13 (i) and (iv)]{LN}, $e(u)$ is degenerate if and only if $\Delta_{\F}$ decomposes with respect to the abelian extension $F(M_\omega)^{G_{34}}\cong F(M_\omega)^{G_{134}}\otimes F(M_\omega)^{G_{234}}$.
\endproof

\begin{remark}  Setting $N_2=F(M_\omega)^{G_{34}}$ and combining Corollaries \ref{cor2} and \ref{cor3} we see that there exist two maximal abelian subfields $N_1$ and $N_2$ of $\Delta_{\F}$ such that 
\[\Delta_{\F} \in \dec(N_1/F(M_\omega)^G)\] and 
\[\Delta_{\F} \notin \dec(N_2/F(M_\omega)^G).\]
\label{remdec}
\end{remark}

\begin{remark}Let $\Delta=(K/F, G, v,d)$ be an abelian crossed product defined by the finite abelian group $G$ of rank $r$, the matrix $v \in M_r(K^*)$ and the vector $d=(d_i)\in (K^*)^r$.  As mentioned in the introduction the generic abelian crossed product associated to $\Delta$ is the abelian crossed product
\[\mathcal{A}_{\Delta}=(K(x_1,\ldots,x_r)/F(x_1,\ldots,x_r),G,v,dx)\]
where the $x_i$ are independent indeterminates and $dx=(d_ix_i)_{i=1}^r\in(K(x_1\ldots,x_r)^*)^r$.  $\mathcal{A}_{\Delta}$ is a division algebra of index $|G|$ and exponent $\textrm{lcm}(\exp(G),\exp(\Delta))$.
Independently in  \cite[Theorem 2.3.1]{mck-galoissubfields} (in the case $\chara(F)=p$) and \cite[Theorem 3.5]{mounirh1} $\mathcal{A}_{\Delta}$ is shown to be indecomposable if $\Delta$ is defined by a non-degenerate matrix.  Therefore, for $\Delta_{\F}$ as in Theorem \ref{p1}, the generic abelian crossed product $\mathcal{A}_{\Delta_{\F}}$ is indecomposable of index $p^2$ and exponent $p$.  This example is in contrast to the indecomposable generic abelian crossed product examples given in \cite[section 3.3]{mck-galoissubfields}.  In those examples the generic abelian crossed products $\mathcal{A}_{\Delta}$  defined by $\Delta=(K(x_1,\ldots,x_r)/F(x_1,\ldots,x_r),G,v,dx)$ have the property that $\Delta$ itself is indecomposable (\cite[Remark 3.3.3]{mck-galoissubfields}).  Moreover in that case the defining matrix $v$ is shown to be non-degenerate by considering the torsion in $\textrm{CH}^2(SB(\Delta))$, the chow group of co-dimension 2 cycles of the Severi-Brauer variety of $\Delta$ and applying work of Karpenko \cite[Prop. 5.3]{Karpenko2}.  Though the calculations in the proof of Theorem \ref{p1} may have been tedious, they are elementary is nature, and do not appeal to the torsion in $\textrm{CH}^2(SB(\Delta_{\mathcal{F}_A}))$.
\end{remark}

\begin{remark} Since the abelian crossed product $\Delta_{\F}$ is decomposable and defined by a non-degenerate matrix, decomposability of abelian crossed products is not determined by degeneracy of the matrix defining them (see \cite[Prop. 3.1.1]{mck-galoissubfields} and remark \ref{remark3}). 
\end{remark}

\bibliographystyle{alpha} 
\bibliography{mckinniebib}
\end{document}